\documentclass[12pt]{article}
\usepackage{amsmath,amssymb}
\newcommand{\G}{\Gamma}
\newcommand{\nin}{\noindent}
\newcommand{\vs}{\vspace*}

\newcommand{\eset}{\emptyset}

\newcommand{\seq}{\subseteq}

\newcommand{\llrarrow}{\longleftrightarrow}
\newcommand{\lrarrow}{\leftrightarrow}

\newcommand{\ol}{\overline}
\newcommand{\ul}{\underline}

\title{\bf SIMPLE GRAPHS AND COMMUTATIVE ZERO-DIVISOR SEMIGROUPS}
\author{Tongsuo Wu\thanks{Corresponding author. {\small email: tswu@sjtu.edu.cn}}
 and Li Chen\thanks{chenli830202@sjtu.edu.cn}\\
Department of Mathematics\\  Shanghai Jiaotong University\\
Shanghai 200030, P. R. China\\}
\date{}
\begin{document}
\baselineskip=16pt \maketitle \vspace{4mm}
\begin{center}
\begin{minipage}{12cm}

\vs{3mm}\begin{center} {\small\bf ABSTRACT}\end{center}

\vs{3mm}\nin {\small In this paper, we study commutative
zero-divisor semigroups determined by graphs. We prove a uniqueness
theorem for a class of graphs. We show two classes of graphs that
have no corresponding semigroups. In particular, any complete graph
$K_n$ together with more than three end vertices and any complete
bipartite graph together with more than one end vertices have no
corresponding semigroups. We also determine all possible
zero-divisor semigroups whose zero-divisor graph is the complete
graph $K_3$ together with two end vertices. }

\vs{3mm}\noindent {\small {\it Key Words:} Commutative zero-divisor
semigroup, Simple connected graphs }

\vs{3mm}\nin {\small\it 1991 Mathematics Subject Classification:}
{\small 20M14, 05C90}

\end{minipage}
\end{center}

\vs{3mm} \begin{center}1.\, {\small INTRODUCTION}\end{center}

\vs{3mm} For any commutative semigroup $S$ with zero element $0$,
there is an undirected zero-divisor graph $\G(S)$ associated with
$S$ ([8]). The vertex set of $\G(S)$ is the set of all nonzero
zero-divisors of $S$,  and for distinct vertices $x$ and $y$ of
$\G(S)$, there is an edge connecting $x$ and $y$ if and only if
$xy=0$. In DeMeyer et al (2005, 2002),  some fundamental properties
and possible structures of $\G(S)$ were studied. For example, for
any semigroup $S$, it was proved that $\G(S)$ is a connected simple
graph with diameter less than or equal to 3, and that the core of
$\G(S)$ is a union of triangles and squares while any vertex of
$\G(S)$ is either an end vertex or in the core, {\it if there exists
a cycle in $\G(S)$}. Many examples of graphs were given in DeMeyer
et al (2005, 2002) and Zuo et al. (2005) to give positive or
negative answers to the following general problem: Given a connected
simple graph $G$, does there exist a semigroup $S$ such that
$\G(S)\cong G$? The zero-divisor graphs were first studied for
commutative rings, see, e.g., Anderson D. et al.(1993), Anderson
D.F. et al.(1999,2003), Beck (1988), Chen (2003), DeMeyer et
al.(2002), Lu et al.(2004), Mulay (2002) and Wu (2005).

For any semigroup $S$, let $T=Z(S)$ be the set of all zero-divisors
of $S$. Then $T$ is an ideal of $S$ and in particular, it is also a
semigroup with the property that all of its elements are
zero-divisors of the semigroup $T$. We call such semigroups {\it
zero-divisor semigroups}. Obviously we have $\G(S)\cong \G(T)$. For
a given connected simple graph $G$, if there exists a zero-divisor
semigroup $S$ such that $\G(S)\cong G$, then we say that {\it $G$
has corresponding semigroups}, and we call $S$ a {\it semigroup
determined by the graph $G$}. In this paper, we study semigroups
determined by some graph $G$. We first give a class of graphs $\G_n$
such that $\G_n$ has a unique corresponding semigroup for each $n$.
(In Wu (2005, Proposition 3.1) a similar uniqueness result was also
obtained for the directed zero-divisor graphs of noncommutative
rings, and the result plays a key role in Wu (2005).) The previous
works in DeMeyer et al (2005, 2002) and Zuo et al.(2005), as well as
works in this paper show that most possible graphs have multiple
corresponding semigroups. The number of semigroups corresponding to
a graph increases rapidly if one end vertex is deleted. On the other
hand, for a graph $G$ having corresponding semigroups (e.g., the
complete graph $K_n$ together with an end vertex, or the complete
bipartite graph $K_{m,n}$ together with an end vertex), if we add
more than two end vertices to $G$, then the resulting graph may have
no corresponding semigroups, as will be shown in the third section
of this paper. This shows that the correspondence between semigroups
and the possible graphs is rather sensitive.

All semigroups in this paper are multiplicatively commutative
zero-divisor semigroups with zero element $0$, where $0x=0$ for all
$x\in S$, and all graphs in this paper are undirected simple and
connected. For any vertices $x,y$ in a graph $G$, if $x$ and $y$ are
adjacent, we denote it as $x-y$ or occasionally, $x\lrarrow y$. For
other graph notations adopted in this paper, please refer to Wilson
(1985).

\vs{4mm}\begin{center}2.\, {\small COMMUTATIVE SEMIGROUPS DETERMINED
BY SOME SIMPLE GRAPHS}  \end{center}

We begin with the following example.

\vs{3mm}\noindent{\bf Example 2.1} For any finite set $A$ with $n$
elements, say, $A=\{1,2,\cdots, n\}$, there is an associated
commutative semigroup $P_n=\{x_B\,|\,B\subseteq A\}$. The
multiplication of $P_n$ is defined by $x_Cx_B=x_{C\cap B}$. It is
straight forward to verify that $P_n$ is a commutative semigroup
with the identity element $x_A$. Also $x_\emptyset $ is the zero
element of $P_n$, i.e., $x_\emptyset x_B=x_\emptyset$, for each
element $x_B\in P_n$. Then we have semigroup isomorphisms
$P_n\cong (2^A,\cap)\cong (2^A, \cup)$, where $2^A$ is the power
set of $A$. Let $\Bbb Z_2$ be the ring of integers modulo $2$. Let
$\Bbb Z_2^{(n)}$ be the ring direct sum of $n$ copies of $\Bbb
Z_2$, and consider its multiplicative semigroup
$(Z_2^{(n)},\cdot)$. It is easy to verify that the map
$$\sigma: (2^A,\cap)\to (Z_2^{(n)},\cdot), B\mapsto (y_1,y_2,\cdots, y_n),
\textrm{where $y_i$}=\left\{\begin{array}{ll} 0&\textrm{if $i\not\in B$}\\
1&\textrm{if $i\in B$}
\end{array}\right.$$
is a semigroup isomorphism.

Denote by $\G_n$ the zero-divisor graph of $P_n$. $\G_n$ is a
symmetric graph with a moderate number of edges. Below we list
some properties of this graph:

1. $V(\G_n)=P_n-\{0,1\}$ and hence it contains $|V(\G_n)|=2^{n}-2$
vertices.

2. For any $x_B\in \G_n$ with $|B|=i$, let $N(x_B)$ be the
neighborhood of $x_B$, i.e., $N(x_B)=\{y\in \G_n\,|\, \textrm{$x -
y$ in $\G_n$}\}$. Then
$|N(x_B)|=C_{n-i}^1+C_{n-i}^2+\cdots+C_{n-i}^{n-i}=2^{n-i}-1$.

3. The edge number is
$|E(\G_n)|=C_n^12^{n-2}+C_n^22^{n-3}+\cdots+C_n^{n-1}-2^{n-1}+1$.

4. The clique number of $\G_n$ is $n$. When $n\ge 3$, the diameter
of $\G_n$ is $3$ and $\G_n$ has $n$ end vertices.

5. The automorphism group of $\G_n$ is the symmetric group $S_n$.
Thus this graph is highly symmetric.

Especially, $\G_2$ is just the complete graph $K_2$, and $\G_3$ is
the complete graph $K_3$ together with three end vertices linked
with distinct vertices of $K_3$. The graph $\G_4$ has 14 vertices
and $25$ edges. $|V(\G_5)|=30, |E(\G_5)|=90$.

\vs{3mm}For this graph $\G_n$ with a moderate edge set, we have
the following result:

\vs{3mm}\noindent{\bf Theorem 2.2} {\it Let $A=\{1,2,\cdots,n\}$ and
assume $n\ge 3$. If $S$ is a commutative zero-divisor semigroup
whose graph $\G(S)$ is isomorphic to $\G_n$, then $S$ is isomorphic
to the  zero-divisor semigroup $P_n-\{1\}$.}

\vs{3mm}\nin{\it Proof.} Assume that $S$ is a zero-divisor semigroup
such that $\G(S)= \G_n$. By the property of $\G_n$, we can have a
labeling $x$ to the elements of $S^*=S-\{0\}$
$$S=\{x_B\,|\, \textrm{$B$ is a proper  subset of \{$1,2,\cdots, n$\}}.\} $$
such that $0=x_\emptyset$, and that for any distinct elements
$x_B,x_C\in S$, $x_Bx_C=0$ if and only if $B\cap C=\emptyset$,
where $B\not= C$.

(1) For any $d\in A$ and any $D\subset A$ with $\{d\}\subset D$, we
can assume $x_dx_D=x_E$ for some non-empty $E$. If $d\not\in E $, we
have distinct numbers $r,s\not= d$ in $A$  such that $r\in E$. If
$E=\{r\}$, then let $H=\{r,s\}$. If $|E|\ge 2$, then take $H=\{r\}$.
In each case, $E\not= H,E\cap H=\{r\}$, and hence $x_Ex_H\not=0.$ On
the other hand, we have $0=x_Hx_dx_D=x_Hx_E$. This contradiction
shows that whenever $\{d\}\subset D\subset A$, we have an $E$ such
that $x_dx_D=x_E$ and $d\in E$. If there exists some $m\in E-\{d\}$,
then $0\not= x_Ex_m=x_dx_mx_D=0x_D=0$, a contradiction. Thus we
obtain $x_dx_D=x_d$.

(2) For any $d\in A$, now we proceed to show that $x_d^2=x_d$. First
we claim that $x_{dr}x_{ds}=x_d $ holds for distinct elements $r,s,
d\in A$. In fact we can assume $x_{dr}x_{ds}=x_G$ for some non-empty
proper subset $G$ of $A$. If $d\not\in G$, then
$0=x_dx_G=(x_dx_{dr})x_{ds}=x_dx_{ds}=x_d$. So we must have $d\in
G$. If $\{d\}\subset G$, then we take an $m\in G-\{d\}$ and we have
$x_m=x_Gx_m=x_{dr}(x_{ds}x_m)$. Since $r\not=s$, either $r\not=m$ or
$s\not=m$. Either case gives $x_{dr}x_{ds}x_m=0$, a contradiction.
Thus $x_{dr}x_{ds}=x_d $.

Now by taking distinct numbers $r,s\not= d$ in $A$, we have the
equality of $x_d=x_dx_{ds}=(x_dx_{dr})x_{ds}=x_d^2. $ By (1),
$x_dx_D=x_d, \forall d\in D\subset A$.

(3) In the following we want to prove that $x_Bx_C=x_{B\cap C}$
holds for any proper subsets $B,C$ of $A$, and this will prove the
uniqueness of the zero-divisor semigroups $S$ with $\G(S)\cong
\G_n$.

If $B\cap C=\eset$, then we are done. Assume $B\not=C$ and $B\cap
C\not=\eset.$ Then in $x_Bx_C=x_E$, $E\not=\eset$. If $E\not\seq
B\cap C$, there is an element $d\in E-B\cap C$, and assume further
that $d\not\in B$. Then $x_d=x_dx_E=(x_dx_B)x_C=0$, a
contradiction. If $B\cap C\not\seq E$, there exists an element
$d\in B\cap C-E$. Then $0=x_dx_E=(x_dx_B)x_C=x_d$, another
contradiction.

If $B=C\not=\eset$, we obtain
$x_d(x_Bx_C)x_d=(x_dx_B)(x_Cx_d)=x_d^2=x_d$ for any $d\in B$, and
therefore $x_B^2\not=0$. Assume $x_B^2=x_E$, where $E\not=\eset$.
If $B\not\seq E$, then for any $d\in B-E$, we have
$0=x_dx_E=x_dx_Bx_B=x_d$, a contradiction. If $E\not\seq B$, then
for any $r\in E-B$, we have $0=x_rx_B^2=x_rx_E=x_r$, another
contradiction. This shows that $x_B^2=x_B$ for any proper subset
$B$ of $A$, and this also completes the proof.\hfill $\Box$

\vs{3mm}\nin{\bf Remarks} (1) Theorem 2.2 does not hold in the case
of $n=2$. Actually, in Zuo et al. (2005, Proposition 3.1) a
nontrivial commutative zero-divisor semigroup $S$ was constructed
such that $\G(S)\cong K_n$. Thus the above uniqueness results fails
for the complete graphs $K_n$. The reason for this lies in the fact
that $K_n$ has too many edges.

(2) In Wu (2005, Proposition 3.1), a uniqueness result was also
obtained for the directed zero-divisor graphs of noncommutative
rings $R$. The following results were proved: For any ring $R$, if
$\G(R)$ has a source vertex(sink vertex, respectively) $x$ such that
$x^2=0$, then $R$ is uniquely determined. This condition was proved
to be equivalent to the condition that the graph $\G(R)$ has exactly
one source vertex (sink vertex, respectively).

\vs{3mm}Let us now consider all subgraphs of $\G_3$ containing
$K_3$ as a subgraph. There are only four such graphs:

(1) The complete graph $M_{3,0}$, i.e., $K_3=\{a_1,a_2,a_3\}$.
This graph has multiple pairwise non-isomorphic corresponding
commutative zero-divisor semigroups. ([6, 12])

(2) The graph $M_{3,1}$, i.e., the complete graph $K_3$ together
with an end vertex, say, $M_{3,1}=\{a_1,a_2,a_3\}\cup \{x_1\}$,
where $a_1 - x_1$. By [6], the graph $M_{n,1}$ is the graph of a
semigroup for any $n\ge 1$.

(3) The graph $M_{3,2}$: $K_3$ together with two end vertices,
say, $M_{3,2}=\{a_1,a_2,a_3\}\cup \{x_1,x_2\}$, where $a_1 - x_1$
and $a_2-x_2$.

(4) $M_{3,3}=\{a_1,a_2,a_3\}\cup \{x_1,x_2,x_3\}$, where $a_1 -
x_1$,$a_2-x_2$ and $a_3-x_3$. By the previous Theorem 2.2,
$M_{3.3}=\G_3$ and it has a unique corresponding commutative
zero-divisor semigroup.

\vs{3mm} \nin{\bf Theorem 2.3} (1) {\it Each of the four subgraphs
of $\G_3$ containing $K_3$ is the graph of a semigroup.}

(2) {\it $M_{3,2}$ has three pairwise non-isomorphic corresponding
commutative zero-divisor semigroups. All the possible
multiplication tables on $V(M_{3,2})$ are listed in the
following}:

\vs{2mm} $\begin{array}{c|ccccc}
\cdot& a_1&a_2&a_3&x_1&x_2\\
\hline
a_1&0&0&0&0&a_1\\
a_2&&0&0&a_2&0\\
a_3&&&a_3&a_3&a_3\\
x_1&&&&x_1&a_3\\
x_2&&&&&x_2\\
&&\tiny{\text{  Table 2.1}}&
\end{array},\quad $
$\begin{array}{c|ccccc}
\cdot& a_1&a_2&a_3&x_1&x_2\\
\hline
a_1&\text{$0$ or $a_1$}&0&0&0&a_1\\
a_2&&a_2&0&a_2&0\\
a_3&&&a_3&a_3&a_3\\
x_1&&&&x_1&a_3\\
x_2&&&&&x_2\\
&&\tiny{\text{Table 2.2}}& \tiny{\text{and}}&\tiny{\text{Table
2.3}}\end{array}.$

\vs{3mm} \nin{\it Proof.} We need only to prove (2). Assume that
$S=\{0,a_1,a_2,a_3,x_1,x_2\}$ is a commutative semigroup with
zero-element $0$ such that $\G(S)=M_{3,2}$

First we claim that $a_1x_2=a_1$ and $a_2x_1=a_2$. In fact, both
$\{0,a_i\,|\,i=1,2,3\}$ and $\{0,a_j\,|\,j=1,2\}$ are ideals of $S$.
Thus it is easy to show that the equalities hold. In this case, we
also have $x_i^2\not=0, \forall\, i=1,2.$ In a similar way, it is
routine to verify that $\{0, a_i\}$ are ideals of $S$, where
$i=1,2$. Thus $a_i^2=a_i$ or $a_i^2=0$ for $i=1,2$.

Next consider the nonzero element $a_3x_1$. We claim that
$a_3x_1=a_3$. In fact, it is obvious that $a_3x_1\not=x_1,x_2$. Now
assume $a_3x_1=a_1$. Then $a_1^2=0$, $a_3x_1^2=0$, and hence
$x_1^2\not=0,x_1,x_2.$ But $a_2=a_2x_1^2=a_2a_1=0$ if $x_1^2=a_1$,
and $a_2=a_2a_3=0$ if $x_1^2=a_3$. Finally, we must have
$x_1^2=a_2$. In this case, we already know $x_1x_2\not=0$ since
$a_3(x_1x_2)=a_1x_2=a_1$. Then from $0=a_2x_2=x_1^2x_2=(x_1x_2)x_1$,
one obtains $x_1x_2=a_1$ since $x_1x_2\not=0$ and $x_1^2\not=0$.
Finally, one obtains $a_1x_2=a_3(x_1x_2)=a_3a_1=0$, a contradiction.
Thus $a_3x_1\not=a_1$. These discussions show that the product
$a_3x_1$ is either $a_2$ or $a_3$.

For later convenience, we now collect all known facts as follows:
$$a_1x_2=a_1, a_2x_1=a_2, x_i^2\not=0, \textrm{$a_3x_1=a_2$ or $a_3$}.\quad\quad (\triangle)$$

(i) \ul{Case 1}. Assume $a_3x_1=a_2$. Then $a_2^2=0$. Since
$(a_3x_2)x_1=(a_3x_1)x_2=0$ and $x_1^2\not=0$, we obtain ${\bf
a_3x_2=a_1}$ and therefore, ${\bf a_1^2=0}$. Next we claim that
$x_1^2=x_1$. In fact, it is obvious that $x_1^2\not= a_1,a_3$. If
$x_1^2=a_2$, then we obtain $0=a_3x_1x_1=a_2$. If $x_1^2=x_2$,
then $a_1x_2=0$, this is again impossible. Hence ${\bf
x_1^2=x_1}$.

Now consider $x_2^2$. If $x_2^2=a_2$ or $a_3$, then
$a_1=a_1x_2^2=0$. If $x_2^2=a_1$, then $0=a_3x_2x_2=a_1x_2$. If
$x_2^2=x_1$, then $a_1=a_1x_2^2=0$. Thus we obtain ${\bf
x_2^2=x_2}$. Since $x_i^2=x_i$ ($i=1,2$), we have $x_1x_2=a_3$.
Hence $a_3^2=a_3$. But then from $a_3^2x_2=a_3(a_3x_2)=a_3a_1=0$
and $x_2^2\not=0$, we obtain $a_3=a_3^2=a_2$, a contradiction.

In conclusion, under the assumption that $a_3x_1=a_2$, there is no
associative multiplication table. Thus we must have $a_3x_1=a_3$. By
symmetry, we also have $a_3x_2=a_3$.

(ii) \ul{Case 2}. Assume {\bf $a_3x_1=a_3x_2=a_3$}. In this case, we
first observe that $x_1x_2\not=0$ since $a_3x_1x_2=a_3x_2.$ Now that
{\it the equalities in $(\triangle)$ are still valid in this
subcase}, we claim that $x_i^2=x_i$. In fact, since $a_3x_i^2=a_3,$
$a_jx_i^2=a_j$ for $j\not=i$, thus $x_i^2$ is adjacent to neither
$a_3$ nor $a_j$ but it still connects to $a_i$. Thus $x_i^2=x_i$.

Now we want to show that $x_1x_2=a_3$. Since $x_1x_2$ is adjacent to
both $a_1$ and $a_2$, we have $x_1x_2\in \{a_i\,|\,i=1,2,3\}$. If
$x_1x_2=a_1$, then we obtain $x_1x_2=x_1^2x_2=a_1x_1=0$. This shows
that $x_1x_2=a_3$. Hence $a_3^2=a_3$ and thus $\{a_3,0\}$ is an
ideal of $S$.

In summary, the following are the only possible multiplication
tables in Case 2:

\vs{2mm}$\begin{array}{c|ccccc}
\cdot& a_1&a_2&a_3&x_1&x_2\\
\hline
a_1&\text{$0$ or $a_1$}&0&0&0&a_1\\
a_2&&\text{$0$ or $a_2$}&0&a_2&0\\
a_3&&&a_3&a_3&a_3\\
x_1&&&&x_1&a_3\\
x_2&&&&&x_2
\end{array}\quad(\text{Table 2.1, 2.2 and 2.3}).\quad$

\vs{2mm}Since the automorphism group of the graph $M_{3,2}$ has only
two elements, among the four tables, only two tables are isomorphic.
So there are at most three semigroup structure on $S$ such that
$\G(S)\cong M_{3,2}$. The final work is to verify that each table
defines an {\it associative} binary operation on $S$. This is really
the case, by direct verification. This completes the whole
proof.\hfill $\Box$

\vs{3mm}The following Corollary corrects a mistake in Zuo et
al.(2005, Example 2.9):

\vs{3mm}\nin{\bf Corollary 2.4} The following graph $G$ has a
unique corresponding zero-divisor semigroup:
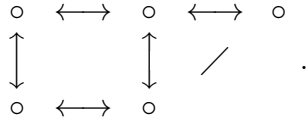
\[\begin{array}{cccccc}
\circ&\llrarrow &\circ&\llrarrow&\circ\\
\Big\updownarrow&&\Big\updownarrow&\diagup\\
\circ&\llrarrow&\circ&&\\
\end{array}.\]

\vs{3mm}\nin{\it Proof.} In case 1 of the proof of Theorem 2.3, we
reverse the procedure in the following way: First deduce
$a_3x_1=a_2,a_3x_2=a_1$. Then verify $a_2x_1=a_2,a_1x_2=a_1$. Then
$x_i^2=x_i$ ($i=1,2$). Finally, $a_i^2=0$. Thus we obtain a unique
associative multiplication table on $S=\{a_1,a_2,a_3,x_1,x_2,0\}$
such that $\G(S)\cong G$. $$\begin{array}{c|ccccc}
\cdot& a_1&a_2&a_3&x_1&x_2\\
\hline
a_1&0&0&0&0&a_1\\
a_2&&0&0&a_2&0\\
a_3&&&0&a_2&a_1\\
x_1&&&&x_1&0\\
x_2&&&&&x_2
\end{array}
\hspace{3cm} \Box$$

\vs{3mm} \nin{\bf Remark.} It is natural to continue the work of
determining all possible semigroup structure for the graph
$M_{3,1}$, which is obtained by deleting one end vertex from the
above graph $M_{3,2}$. While doing so, we have found fifteen
pairwise non-isomorphic associative multiplication tables for the
graph $M_{3,1}$. These tables are too many to be included here.

\vs{4mm}\begin{center}3.\, {\small TWO CLASSES OF GRAPHS WITH NO
CORRESPONDING SEMIGROUPS }\end{center}

\vs{3mm}Consider the following graphs and their generalizations
\[\begin{array}{cccccccccccccccc}
\circ&\llrarrow &\circ&&&&& &&     \circ&\llrarrow &\circ&\llrarrow&\circ&\llrarrow\circ\\
&&\Big\updownarrow&\diagdown&&&&&&                        &     &\Big\updownarrow&\diagup\diagdown&\Big\updownarrow&&\\
\circ&\llrarrow&\circ&\llrarrow&\circ&\llrarrow&\circ&,&&
\circ&\llrarrow&\circ&\llrarrow&\circ&\llrarrow\circ\\
&&&\text{{\tiny Fig.3.1}} &&&&&&&&&\text{{\tiny Fig.3.2}}&&&
\end{array}.\]
The graph in Fig. 3.1 is just the $\G_3$, while the other graph in
Fig.3.2, denoted as $M_4$, is a kind of generalization of $\G_3$.
The core of $M_4$ is the complete graph $K_4$.

By Theorem 2.2, the graph in Fig 3.1 uniquely determines a
zero-divisor semigroup. Now let us consider the graphs $M_n$ for
$n\ge 4$, where $M_n$ is the complete graph $K_n$ together with $n$
end vertices such that each vertex of $K_n$ connects to an end
vertex. Do the graph $M_n$ and it's subgraphs containing $K_n$ have
corresponding semigroups? By DeMeyer et al.(2005,  Theorem 3(1)) and
Zuo et al. (2005, Proposition 3.1), $K_n$ and $K_n$ together with an
end vertex do have corresponding semigroups. But for those having
more than three end vertices, the answer is no. This fact is a
special case of the following general result.

\vs{3mm}\nin{\bf Theorem 3.1} {\it For $n\ge 4$, let a simple
graph $L_n=\{a_1,a_2,a_3,a_4,\cdots, a_n\}\cup (\cup_{i=1}^nX_i)$
be the disjoint union of $n+1$ subsets satisfying the following
conditions:}

(1) {\it $X_i$ $(i=1,2,3,4)$ are nonempty subgraphs}.

(2) {\it Each $a_k$ is adjacent to either $a_3$ or $a_4$, for all
$k=1,2,\cdots, n$.}

(3) {\it $a_i-a_j$ if $i\not=j$ and
$X_i\not=\eset,X_j\not=\eset$.}

(4) {\it $a_j$ is adjacent to each vertex of $X_j$, for any
$j=1,2,\cdots,n$.}

(5) {\it There is no edge linking a vertex in $X_i$ with a vertex
in $X_j$ ($i\not=j$)}.

(6) {\it There is no edge connecting $a_i$ with a vertex in $X_j$
($i\not=j$)}.

(7) {\it The graph structures of the subgraphs $X_i$ can be chosen
freely}.

\nin {\it Then for any commutative semigroup $S$, the zero-divisor
graph $\G(S)$ can not be isomorphic to $L_n$.}

\vs{3mm}\nin{\it Proof.} Suppose that there is a commutative
semigroup $T$ with zero element $0_T$ such that $\G(T)\cong L_n$.
Then $Z(T)=S$ is a commutative zero-divisor semigroup with zero
element $0_T$. Without any loss, we can assume $S=V(L_n)=K_n\cup
(\cup_{i=1}^nX_i)\cup \{0_T\}$.

First we show the following result: If $X_j\not=\emptyset$, then
$\{0,a_j\}$ is an ideal of $S$ and hence,  $a_jx_i=a_i$ for all
$x_i\in X_i$, $j\not=i$. In fact, let $z$ be any vertex. Then
$x_j(a_jz)=0$ for some $x_j\in X_j$. Thus $a_jz\in \{0,a_j\}\cup
X_j$. However there exists some $l$ such that $X_l\not=\emptyset$,
$l\not=j$. Thus $a_jz\not\in X_j$ since by assumption $a_la_jz=0$.
Thus $a_jz\in \{0,a_j\}$ and $\{0,a_j\}$ is an ideal of $S$.

Now we fix two elements $x_1\in X_1,x_2\in X_2$ and consider
$x_1x_2$. Set $z=x_1x_2$. Then $z\neq 0$, and so either $z\in
A=\{a_1,a_2,\cdots,a_n\}$, or $z\in X=\cup_{i=1}^kX_i$. By
assumption (1), we have $a_3z=a_3x_1x_2=a_3\neq 0$ and $
a_4z=a_4x_1x_2=a_4 \neq 0$, and it follows from the assumption (2)
that $z\notin A$. So $z\in X$. But this is contradicting to the fact
of $a_1z=a_1x_1x_2=0$ and $a_2z=0$. \hfill $\Box$

\vs{3mm} The conditions (2) and (3) are unnecessary in Theorem 3.1.
But if we drop one of them, then the result follows immediately from
the known result that the diameter of $\G(S)$ is at most 3..

\vs{3mm}By the result and the proof of Theorem 3.1, we have the
following

\vs{3mm}\nin{\bf Corollary 3.2} {\it The following graphs have no
corresponding semigroups}:

(1) {\it  $M_{n,k}=K_n\cup (\cup_{i=1}^k\{x_i\})$: the complete
graph $K_n=\{a_1,a_2,\cdots, a_n\}$ together with $k$ end vertices
such that $a_i-x_i$ ($n\ge k\ge 4$)}.

(2) {\it Any generalizations of graphs with the following forms
and the {\it related} refinement graphs}:
\[\begin{array}{cccccccccccccccccccc}
&&\circ&&&&&                   &&&\circ&&&&& \\
&&\Big\updownarrow&\diagdown   &&&&&&&\Big\updownarrow&\diagdown\\
\circ&\llrarrow &\circ&\llrarrow&\circ&\llrarrow\circ    &&&\circ&\llrarrow &\circ&\llrarrow&\circ&\llrarrow\circ\\
&&\Big\updownarrow&\diagup\diagdown&\Big\updownarrow&&   &&&&\Big\updownarrow&\diagup\diagdown&\Big\updownarrow&&\\
\circ&\llrarrow&\circ&\llrarrow&\circ&\llrarrow\circ     &&&\circ&\llrarrow&\circ&\llrarrow&\circ&\llrarrow\circ\\
&&&\diagdown&\Big\updownarrow&                           &&&&&&\text{{\tiny Fig.3.4}}&&&\\
&&&&\circ& \\
&&&\text{{\tiny Fig.3.3}}&&&
\end{array}.\]

\vs{3mm} We remark that for each $n\ge 4$, {\it the graph
$M_{n,2}$ has a unique zero-divisor semigroup, while $M_{n,3}$ has
no corresponding zero-divisor semigroups.} These results and the
related constructions and proofs will appear in a subsequent
paper.

\vs{3mm}By DeMeyer et al.(2005,  Theorem 3(2)), any complete
bipartite graph and any complete bipartite graph together with an
end vertex is the graph of a semigroup. The complete bipartite graph
case was also independently discovered in Zuo et al.(2005,
Proposition 3.2). Like the complete graph case, a step further can
lead to a negative result.

\vs{3mm}\nin{\bf Theorem 3.3}  {\it For any $m,n\ge 2$, let the
connected simple graph
$$H_{m,n}=\{a_1,\cdots,a_m;b_1,\cdots,b_n\}\cup X_1\cup Y_1$$ be the
disjoint union of three non-empty subsets satisfying the following
conditions}:

(1) {\it $a_1-b_k$, $a_2- b_k$, $b_1- a_r$ and $b_2- a_r$ for all
possible $k$ and $r$}.

(2) {\it $a_1$ ( $b_1$ ) is adjacent to each vertex of $X_1$ (
$Y_1$, respectively)}.

(3) {\it There is no edge linking a vertex in $X_1$ with a vertex
in $Y_1$}.

(4) {\it There is no edge connecting $a_i$ ( $b_i$ ) with a vertex
in $Y_1$ ( $X_1$, respectively)}.

(5) {\it There is no edge linking $a_i$ ( $b_i$ ) with a vertex in
$X_1$ ($Y_1$, respectively ), for all $i\not=1$}.

(6) {\it $\{a_1,\cdots,a_m;b_1,\cdots,b_n\}$ is a subgraph of the
complete bipartite graph $K_{m,n}$.}

(7) {\it The graph structures of the subgraphs $X_1$ and $Y_1$ can
be chosen freely}.

\nin {\it Then for any commutative semigroup $S$, the zero-divisor
graph $\G(S)$ can not be isomorphic to $H$}.

\vs{3mm}\nin{\it Proof.} Assume in the contrary that there exists a
commutative zero-divisor semigroup $S$ such that $\G(S)\cong
H_{m,n}$. First we conclude that $a_1y=a_1$ and $b_1x=b_1$ for all
$x\in X_1, y\in Y_1$. In fact it is riutine to verify that both
$\{0,a_1\}$ and $\{0,b_1\}$ are ideals of $S$. Especially, we
further deduce that $x^2\not=0$ and $y^2\not=0$ for all $x\in X_1,
y\in Y_1$.

Now we fix two elements $x\in X_1$ and $y\in Y_1$, and we consider
$xy$. If $xy=a_i$ ($i\not=1$), then $a_1a_i=0$, a contradiction.
Similarly, we also have $xy\not= b_j$ for any $j\not=1$. If
$xy=x_1\in X_1$, then $b_1x_1=0$, a contradiction. Similarly, we
also have $xy\not\in Y_1$. This proves that either $xy=a_1$ or
$xy=b_1$. By symmetry, we can assume {\bf $xy=a_1$}. Then we
obtain $a_1^2=0$, and $0=a_1x=x^2y$. Since $y^2\not=0$, we have
either $x^2=b_1$ or $x^2=y_1\in Y_1$. But if $x^2=y_1$, then
$a_1y_1=a_1x^2=0$,a contradiction. Thus $x^2=b_1$ and thus
$b_1^2=b_1x^2=b_1$. From $xy=a_1$, we also obtain
$0=b_ra_1=(b_ry)x$ ($r\ge 2$). Thus we have $b_ry=z_r$ where
either $z_r=a_1$ or $z_r=x_2\in X_1$.

Finally, we consider $a_2y$. If $a_2y=a_r$ ($r\ge 1$), then by
assumption $0=b_2a_r=(b_2y)a_2=z_2a_2$, where either $z_2=a_1$ or
$z_2=x_2\in X_1$, a contradiction in either case. If $a_2y=b_1$,
then $b_1^2=0$, contradicting to the result of $b_1^2=b_1$. If
$a_2y=b_j$ for some $j\ge 2$, then $b_1b_j=0$, contradicting to
the assumption. Finally, we must have $a_2y\in X_1\cup Y_1$. But
if $a_2y=x_1\in X_1$, then $b_1x_1=0$. If $a_2y=y_1\in Y_1$, then
$b_2y_1=0$. In conclusion, $a_2y\not\in S$ and this also completes
our proof. \hfill $\Box$

\vs{3mm}\nin{\bf Corollary 3.4} {\it For any $m\ge 2$ and $n\ge
2$, let $L_{m,n}$ be the complete bipartite graph $K_{m,n}$
together with at least two end vertices which connect to distinct
vertices of $K_{m,n}$. Then $L_{m,n}$ has no corresponding
semigroups. Especially, the following graph has no corresponding
semigroups}:
\[\begin{array}{cccccccc}
\circ&\llrarrow &\circ&\llrarrow&\circ\\
\Big\updownarrow&&\Big\updownarrow&&\\
\circ&\llrarrow &\circ&\llrarrow&\circ\\
&&\text{{\tiny Fig.3.5}} &&
\end{array}.\]

\vs{3mm}\nin{\it Proof.} If at least two end vertices, say, $x$ and
$y$, connect to one part of $K_{m,n}$, then the distance from $x$ to
$y$ is $4$. Thus this $L_{m,n}$ has no corresponding semigroups.

The other case is that $L_{m,n}$ has exactly two end vertices and,
these end vertices connect to different parts of $K_{m,n}$. In
this case, the result follows from Theorem 3.3. \hfill $\Box$.

\vs{3mm}We remark that many negative graphs $G$ of Theorem 3.1 and
Theorem 3.3 satisfy the conditions (1) to (4) of DeMeyer et
al.(2005, Theorem 1), i.e., the graphs $G$ satisfy the following
conditions: (1) $G$ is connected and simple, and the diameter of $G$
is at most three. (2) The core of $G$ is a union of triangles and
squares, and any vertex not in the core is an end-vertex. (3) For
any non-adjacent vertices $x,y$, there exists a vertex $z$ such that
$N(x)\cup N(y)\subseteq \ol{N(z)}$.

We believe that Example 2.1 is an important positive graph in the
theory of zero-divisor graphs of semigroups. In it there is a
complete subgraph $K_n$ with $n$ end vertices linked to $n$ vertices
of $K_n$ respectively, for all $n$. We also noticed that some
negative graphs in Theorem 3.1 (e.g., $M_n$) and Theorem 3.3
contains a cycle $C_n$ ($n\ge 4$)£¬at least four vertices of which
linked to one end vertex respectively. These facts show that the
structures of the zero-divisor graphs of semigroups are complicated.

\vs{4mm}\begin{center}{\bf ACKNOWLEDGEMENT}\end{center}

\vs{3mm} The authors express their appreciation to Dr. Lu Dancheng
who suggests a simplified proof to Theorem 3.1.

\vs{6mm}\begin{center}{\bf REFERENCES}\end{center}

\vs{3mm}

\nin Anderson D.D. and Naseer M. (1993). Beck's coloring of a
commutative ring. {\it J.

\hspace{0.5cm}Algebra.} 159:500-514.

\nin Anderson D.F.; Livingston P.S. (1999). The zero-divisor graph
of a commutative

\hspace{0.5cm}ring. {\it J. Algebra.} 217:434-447.

\nin  Anderson D.F.; Ron Levy; Shapiro J. (2003). Zero-divisor
graphs, von Neumann

\hspace{0.5cm}regular rings, and Boolean algebras. {\it J. Pure
Applied Algebra.} 180:221-241.

\nin Beck I. (1988). Coloring of commutative rings.  {\it J.
Algebra.} 116:208-226.

\nin Chen P.W. (2003). A kind of graph structure of rings. {\it
Algebra Colloquium.}

\hspace{0.5cm}10:229-238.

\nin DeMeyer F.R.; DeMeyer L. (2005). Zero divisor graphs of
semigroups. {\it J. Algebra.}

\hspace{0.5cm}283:190-198.

\nin DeMeyer F.R.; Schneider K. (2002). Automorphisms and
zero-divisor graphs of

\hspace{0.5cm}commutative rings. {\it Internat. J. Commutative
Rings.} 1:93-106.

\nin DeMeyer F.R.; McKenzie T.; Schneider K. (2002) The zero-divisor
graph of a

\hspace{0.5cm}commutative semigroup. {\it Semigroup Forum.}
65:206-214.

\nin Lu D.C.; Tong W.T. (2004). The zero-divisor graphs of abelian
regular rings.

\hspace{0.5cm}{\it Northeast. Math. J.} 20:339-348.

\nin Mulay S.B. (2002). Cycles and Symmetries of Zero-divisors. {\it
Comm. Algebra.}

\hspace{0.5cm}30:3533-3558

\nin Wilson R.J. (1985). {\it Introduction to Graph Theory.} Longman
Inc. New York,

\hspace{0.5cm}Third Edition.

\nin Wu T.S. (2005) On directed zero-divisor graphs of finite rings.
{\it Discrete Math.}

\hspace{0.5cm}296:73-86. DOI:10.1016/j.disc.2005.03.006.

\nin Zuo M.; Wu T.S. (2005) A new graph structure of commutative
semigroups.

\hspace{0.5cm}{\it Semigroup Forum.} 70:71---80.  DOI:
10.1007/s00233-004-0139-8

\end{document}